\documentclass[11pt,a4paper]{article}
\usepackage[utf8]{inputenc}
\usepackage[T1]{fontenc}
\usepackage{lmodern}
\usepackage[english]{babel}
\usepackage{microtype}
\usepackage[a4paper,margin=2.8cm]{geometry}
\usepackage[hidelinks]{hyperref}

\title{The Role of Mathematicians in Mathematics}
\author{Carlo Gasbarri}
\date{}

\makeatletter
\newcommand{\address}[1]{\gdef\@address{#1}}
\newcommand{\email}[1]{\gdef\@email{#1}}
\newcommand{\printaddress}{%
  \par\vspace{2em}\noindent
  \begin{minipage}{\linewidth}\small
  \@address\par
  \emph{E-mail address}: \texttt{\@email}
  \end{minipage}}
\makeatother

\address{Carlo Gasbarri, IRMA, UMR 7501\\
 7 rue Ren\'e-Descartes\\
 67084 Strasbourg Cedex}
\email{gasbarri@math.unistra.fr}
\date{\today}

\begin{document}
\maketitle

\vspace{-1em}

\begin{abstract}
\noindent Some reflections on the role of mathematicians in an age of proofs found and verified by machines.
\end{abstract}
\vspace{0.5em}

André Weil is said to have remarked that, as a young man, he hoped to prove the Riemann Hypothesis; later in life he would have been happy simply to understand a proof of it; and in old age he would have settled for knowing that one existed. The anecdote is usually told for its gentle irony about ambition and age, but it also raises a less comfortable question. What, exactly, do we gain when a theorem is proved? Take the most celebrated case of recent decades. Has the fact that the equation $x^n + y^n = z^n$ has no solutions in positive integers for $n > 2$ really changed anyone's life, or anything substantial in it? One might object that it certainly changed the lives of Andrew Wiles, Richard Taylor, and the many mathematicians directly involved in the proof and in its long aftermath, down to Kevin Buzzard, who has set out to have the whole argument verified by a computer. Their careers have brought academic prestige, prizes and honours, and prominent positions. Look more closely, though, and it is not the absence of solutions that transformed their lives. It is the fact that \emph{they} found the proof, that \emph{they} understood why no solutions exist, and that \emph{they} were able to explain it to the rest of the mathematical community.

I belong to a different and far larger group. Like most professional mathematicians (and I mean \emph{most}: the overwhelming majority of us, not some unlucky fringe) I have had neither the strength nor the good fortune to prove a major conjecture. My work, then and now, has been of a humbler kind: pushing a theory a little further, making it evolve a bit, coming to understand it properly, and, quite often, explaining it to other mathematicians. This is not a confession of failure; it is simply the ordinary condition of the profession. Over the years, watching friends and colleagues reach results that will carry their names, I have more than once asked myself whether I really had a role in the mathematical community at all. At some point I understood that I did. My role is to understand, and to help others understand, why so many mathematical results matter. Without people like me, the so-called geniuses and heroes of mathematics would have no competent audience, no one able to truly appreciate what they have discovered.

There have been many attempts to say what is important in mathematics. One familiar answer appeals to the verdict of time: theories developed long ago, often with no thought of use, have turned out to be of pressing relevance today and to have countless applications. The favourite example is almost too good. In \emph{A Mathematician's Apology} (1940), G.~H. Hardy took a quiet pride in the uselessness of number theory, which he held up as the purest and most innocent of disciplines, one that no war could ever put to use. Fate, it seems, has a sense of humour. A few decades later number theory had become the backbone of modern cryptography. Today it guards banks, governments and military communications, and Hardy's harmless pastime turns out to be a matter of national security.

Hardy's story is real, and it is often told. But it is the exception that everyone remembers precisely because it is an exception. For every Hardy whose ``useless'' mathematics later found its way into the world, there are countless theories that remained, and will most likely remain, exactly as useless as their authors cheerfully admitted. For most of mathematics, the argument from posterity is less a justification than a hope. I cannot help thinking that those who rely on it are a little like the obscure poet who consoles himself with the dream that, centuries from now, someone will stumble upon his forgotten manuscripts, and that he will be revealed, posthumously and gloriously, as the greatest poet of his age. It is a comforting thought. It is also one that, by its very nature, nobody alive will ever be in a position to check.

What Hardy's famous remark really reveals lies elsewhere, and the history of physics tells the same story. Max Planck was certainly not thinking of the countless applications that quantum physics would one day have. He was trying to make sense of how a heated body radiates. Einstein, for his part, did not develop relativity so that we might one day find our way with a GPS receiver, even though that device depends on his theory to work at all. What drove these scientists, and Hardy with them, was the desire to understand nature and the mathematics they were studying. Their interest was never merely to know \emph{whether} something was true. It was to understand \emph{why} it was true, and to discover \emph{how} it could be proved.

In the last few months, however, the question has returned in a far more urgent form, because artificial intelligence has burst into mathematics. On 19 July 2026, Levent Alpöge, a mathematician working at Anthropic, announced in a laconic post on X an explicit counterexample to the Jacobian Conjecture in dimension three, found with the help of the AI model Claude Fable~5~[1]. On 8 September, OpenAI announced that one of its unreleased internal models had settled the Navier--Stokes existence and smoothness problem, one of the Clay Millennium Prize Problems, by showing that an initially smooth flow can develop a singularity in finite time~[2]. On 6 October the same company released a repository of 722 manuscripts, grouped into 372 ``families'' of results, which apparently contain the resolution of a remarkable number of open conjectures~[3]. Finally, on 20 November, Alpöge posted on his personal website a four-page note that reads rather like a stream of consciousness. In it he states that, working with Claude, he believes he has found a proof of the Riemann Hypothesis, and indeed of the Generalized Riemann Hypothesis~[4]. The proof has been released only as Lean code. Within hours, independent mathematicians reported that it compiles and rests only on the standard axioms, and a footnote in the note mentions that a written version generated along the way ran past 250{,}000 pages~[5]. Weil, it seems, may yet be granted the wish of his old age: to know that the Riemann Hypothesis has been proved. Whether anyone will ever \emph{understand} the proof is another matter entirely.

It seems that in the last six months the entire mathematical community has passed from the first to the third stage of André Weil's life, and that it has done so almost without noticing it was skipping the second. If further proof of Weil's extraordinary insight into the world were still needed, here it is: he seems to have anticipated not only the arc of a mathematician's life, but the fate of mathematics itself. The first successes of artificial intelligence could still be checked by any one of us. The counterexample to the Jacobian Conjecture fits in a few lines, and an afternoon with pencil and paper, or a few seconds on any computer algebra system, is enough to convince oneself that it works. What is offered to humanity now is something quite different: an oracle, a \emph{vates} in the old sense of the word, which tells us that the Riemann Hypothesis is true. We cannot follow the argument, and we cannot hope to read it. We can only accept the verdict, together with the assurance of another machine that every step has been checked. Weil, in his old age, would have been content to know that a proof existed. But is this really the condition to which we aspire? Does humanity, or at least the mathematical community, now wish to regard itself as an old man in retirement, content to be told that the problems of its youth have been solved by someone else?

Fortunately, part of the mathematical community seems to be reacting to Alpöge's announcement. On the very same day, Ravi Vakil made public a 940-page preprint by Alessandro Masullo, \emph{The Music of the Primes}, which proposes a proof of the Generalized Riemann Hypothesis, the fruit of some twelve years of work and, in Vakil's words, entirely human. Masullo's work has not been free of controversy. For years his alleged proof circulated more as a rumour than as a text, its hurried release has been greeted with open scepticism, and within hours readers were already pointing out possible gaps. But all of this is precisely the point. What is happening around Masullo's manuscript is a human activity, made of debate and understanding. Mathematicians are reading it page by page, raising objections, proposing repairs, arguing over a single inequality, even when they enlist machines to help them. Perhaps Masullo is wrong; perhaps he is not. What matters is that human beings go on understanding and discussing questions that enlarge the knowledge of humankind, and that whatever the outcome, those who took part will come away knowing more than they did before. Consider the alternative. Suppose that one day a program reports that the proof produced with Claude is wrong, or incomplete: a statement formalized incorrectly, a hidden assumption, a step that does not quite say what it claims. What could we do then? Very little. We could only shrug, say ``never mind, there is nothing I can do about it'', and forget the whole affair. Nobody would have learned anything, because nobody had understood anything in the first place. And suppose, on the contrary, that everything holds: that experts agree the formal statement really is the Riemann Hypothesis, and that Lean certifies it as true beyond any reasonable doubt. What will we have learned about the Riemann Hypothesis? That it is true, and nothing more. We will know no more than before about \emph{why} the zeros of the zeta function lie on the critical line, what hidden structure forces them there, or what this tells us about the primes. The most famous question of our discipline will have been answered, and we will be exactly as ignorant as we were the day before.

The history of mathematics is full of mistakes that turned out to be more fruitful than many correct theorems. In 1848 Jakob Steiner asked how many conics are tangent to five given conics in general position, and answered 7{,}776. The answer was wrong. Sixteen years later Michel Chasles showed that the correct number is 3{,}264, and the effort to understand why Steiner had gone astray, and how such numbers should be counted at all, stands near the very beginning of what would become intersection theory. The episode is so emblematic that David Eisenbud and Joe Harris chose that number as the title of their book, \emph{3264 and All That: A Second Course in Algebraic Geometry}. A century later, in the 1950s, the Belgian mathematician Lucien Derwidué published a proof of the resolution of singularities of algebraic varieties in arbitrary dimension. The proof was flawed. Beniamino Segre and Oscar Zariski argued at length against it, and Zariski dissected it point by point in a review that remains a small masterpiece of the genre. It is hard to believe that this close engagement with a wrong proof did not help them, and through them the whole community, to understand the problem better, a decade before Heisuke Hironaka finally settled it in characteristic zero. And then there is Fermat himself, with whom this essay began. Had he not scribbled in the margin of his copy of Diophantus that he had a truly marvellous proof which the margin was too narrow to contain, a great deal of mathematics, from Kummer's ideal numbers to the modularity of elliptic curves, might never have been developed. Rarely has a claim that was almost certainly mistaken been so productive. Go further back still, and we find the search, lasting nearly two thousand years, for a ``proof'' of Euclid's fifth postulate. It was a long sequence of errors, near misses and ingenious failures, from Proclus to Saccheri and Lambert, and from it eventually came the non-Euclidean geometries of Gauss, Bolyai and Lobachevsky, and with them modern geometry. Or consider the quadrature of the circle, a problem littered for centuries with mistakes, many of them by amateurs and some of them by eminent men. It ended only in 1882, when Lindemann proved that $\pi$ is transcendental, and in doing so gave a decisive impulse to the theory of transcendental numbers. In every one of these cases, the error was not a dead end. It was the beginning of a conversation, and it was the conversation, carried on by human beings who wanted to understand, that made mathematics grow. Not one of these errors would have taught us anything had it been reported to us by a machine as simply \emph{false}.

None of this is meant as an apology for error, still less as a celebration of human weakness. The point is a different one. Mathematical knowledge is not an encyclopaedic list of true statements, to which each new theorem adds one more line. It is the progress of a debate, the slow and collective acquisition of understanding, in which a wrong proof, an ambiguous note in a margin or a centuries-long failure may matter as much as a correct result, because each of them forces us to understand something we had not understood before. J.~Robert Oppenheimer described this movement beautifully in his essay ``Analogy in Science.'' We begin, he observed, with an intuition. We then study it, work on it and make it rational, and in the process the intuition loses much of its original charm. But we gain in exchange a higher intuition, from which our growth in knowledge can begin again. This, it seems to me, is what mathematics truly is: not the possession of truths, but the endless passage from one intuition to a deeper one. And this is precisely what an oracle, however infallible, cannot give us. A statement certified as true, but understood by no one, does not raise our intuition by a single degree. It closes a question without opening a new one.

This note is not meant as yet another indictment of the use of artificial intelligence in mathematics. AI is here, and it is a new tool, a remarkably useful and powerful one. We have seen tools change the way mathematics is done before. Older colleagues still remember the hours spent in the library copying out by hand the latest articles as they arrived. They will tell you, with some nostalgia, that they learned a great deal that way, unlike the young, who do nothing but download papers they will never read. Further back, Euler, blind in the last part of his life, carried out long and intricate computations in his head and dictated them to his assistants, and his productivity seems to have suffered little from it. These practices, and many others like them, have been swept away by the changing times, and the changing times have now brought us to ``conversations'' with AI. Such conversations already make it possible to take interesting and often unexpected steps forward in our research. For those of us whose role, as I have described it, is to understand and to explain, they can be precious in another way. They help us make sense of difficult results, fill in the gaps that authors leave for the reader, and check long and delicate computations. What I want to resist is something else entirely. It is the growing tendency to leave human beings out of the evolution of knowledge, as if humanity's part in the advancement of culture and understanding were a marginal one, that of a spectator watching from the stands while the game is played by others.

The scientists and mathematicians who have built modern language models and proof assistants are doing extraordinary work, and it is already changing almost every field of human activity; I say this with genuine admiration. What is at stake is keeping the human being at the centre of science, of mathematics and of reason itself. For some five centuries, let us say since Galileo, science has done exactly that. It left behind the medieval \emph{ipse dixit}, the principle that a statement is true because Aristotle, or some other authority, has said so, and put human reason and human understanding in its place. A proposition was to be accepted not because of who uttered it, but because anyone willing to look through the telescope, repeat the experiment or follow the proof could see for themselves that it held. It would be a strange irony if, after five hundred years, we returned to an \emph{ipse dixit} of a new kind, in which the authority is no longer a philosopher but a machine. Let me be clear: I believe we are entitled to use whatever method we like to find the proof of a theorem or to build a theory. It may come at a desk, after months of stubborn effort, alone or as part of a team. It may come in the middle of heavy physical labour: Chebotarev is said to have found the key to his celebrated density theorem not in a study but with his arms full, during a move, which suggests that inspiration does not always wait for the mind to be at rest. It may come in a dream, which seems to have happened to more mathematicians than one might suppose. It may come from a dialogue with a machine (\dots). How a result is found is nobody's business but its discoverer's. What I consider essential is something else. A result, an article or a theory must have one or more human beings who answer for it: people who can explain it, defend it against objections, and take it further. The counterexample to the Jacobian Conjecture is a good illustration. Found with the help of a machine, it did not remain a mysterious verdict for long. Within days mathematicians had extended it to an infinite family and given it a geometric explanation~[7], and Arno van den Essen, who had worked on the conjecture since 1985, wrote a note to understand how it could have been found without AI at all~[6]. The result was digested, made human, so to speak, and in this way it has become part of human knowledge. Moritz Schlick, echoing Wittgenstein, held that the meaning of a proposition is the method of its verification. If so, a theorem verified only by a machine, through a process that no human being can follow, has a meaning for that machine, perhaps, but not yet for us. And there is more at stake than the meaning of a single theorem. Mathematics has always moved, and must go on moving, towards simplification. What was once difficult, the hard-won achievement of a handful of experts, becomes with time a standard chapter, then a classroom topic, and sometimes no more than an exercise at the end of a textbook. Results that once took a lifetime to reach are now proved by students in an afternoon. This is not a loss of depth. It is the sign that mathematical ideas have been truly assimilated by humankind, digested and made part of the common heritage of anyone who cares to learn them. It is essential that this process continue. We must not accept that machines bring it to a halt simply by handing us guarantees of truth, however irrefutable those guarantees may be. A proof of 250{,}000 pages that no one will ever simplify, because no one will ever understand it, is not the end point of this long human process. It is its interruption.

We human beings want, and must, keep our central place in knowledge. We must demand, and accept as a principle, that a theorem, a theory or a result, however it was found, becomes part of our knowledge only when we are able to \emph{understand} it, \emph{explain} it, \emph{develop} it and \emph{digest} it. That understanding may come through discussion, sometimes heated, through conversation, through explanation, through any of the countless ways in which mathematicians have always argued with one another. It may well come with the help of a machine. But in no way should it be imposed on us by an oracle. It has often been said that mathematics is a profoundly democratic discipline. Its doors are open to anyone willing to enter, to study and to make the effort to understand, whatever their origin, their means or their titles. No authority, however eminent, can make a statement true. It must be shown to be true, in a way that anyone who takes the trouble can follow. A result that no one can understand, certified by a voice that no one can question, would be the first true exception to this principle in the history of our discipline. It is our duty as human beings to go on demanding that mathematics remain open to all who wish to understand it, and to refuse to let it become the private conversation of machines.

Which brings me, at last, to the title of this essay: the role of mathematicians in mathematics. I hope I have made my position clear. Mathematicians play in mathematics the same role that human beings play in humanity. The statement sounds like a tautology, and in a sense it is one, which is exactly the point. Without mathematicians, mathematics would not merely lose its meaning; it would not exist at all. Physicists like to say that it makes no sense to ask what came before the Big Bang, because time itself began with it. Scientists understand the argument perfectly well. The man in the street does not, and goes on asking, with stubborn common sense, what there was \emph{before}. Today, oddly enough, the man in the street has acquired a new and rather distinguished company. He now includes the most ardent enthusiasts of artificial intelligence, often brilliant people with impressive credentials, who look forward to a future in which proofs will be produced by machines, checked by other machines, and certified as true without any human being ever needing to understand them. They imagine a mathematics that would go on growing, theorem after theorem, in a world where no one is left to understand it. Like the man in the street asking what came before the Big Bang, they are asking a question that has no meaning. A theorem that nobody understands, proved by one machine and verified by another, is no more a part of mathematics than a book written in a language no one will ever speak is a part of literature. And the more oracles we have, the more we may need people able to turn their verdicts into understanding, and to share that understanding with others. Mathematics without mathematicians is not a more efficient mathematics. It is simply nothing at all.

\section*{References}
\begin{enumerate}
\renewcommand{\labelenumi}{[\arabic{enumi}]}
\item L.~Alpöge, post on X (@\_\_alpoge\_\_), July 2026.
\item OpenAI, ``On the Navier--Stokes Millennium Prize Problem,'' \url{openai.com}, 8 September 2026.
\item OpenAI, mathematics repository, \url{github.com/openai/math}, released 6 October 2026.
\item L.~Alpöge, note on his personal website, 20 November 2026.
\item ``Anthropic mathematician claims to prove Riemann hypothesis, most famous problem in all of math,'' \emph{Scientific American}, 20 November 2026.
\item A.~van den Essen, ``An elementary way to find a counterexample to the Jacobian Conjecture,'' arXiv:2609.17795, September 2026.
\item S.~Gao, ``Counterexamples to the Jacobian conjecture in dimensions greater than two,'' arXiv:2608.00222, July 2026 (revised November 2026).
\end{enumerate}

\printaddress

\end{document}